\documentclass[a4paper, 11pt]{article}
\usepackage[utf8]{inputenc}
\usepackage{fullpage}
\usepackage[all]{xy}
\usepackage{color}
\usepackage[english]{babel}
\usepackage[numbers]{natbib}

\usepackage{hyperref}
\usepackage{amsmath}
\usepackage{amssymb}
\usepackage{amsthm}
\usepackage{eucal}
\usepackage{mathrsfs}
\usepackage{mathtools}
\usepackage[disable]{todonotes} 

\renewcommand{\phi}{\varphi}

\newcommand{\frk}[1]{\mathfrak{#1}}

\newcommand{\Z}{\mathbb{Z}}

\newcommand{\C}{\mathbb{C}}
\newcommand{\R}{\mathbb{R}}

    {\@beginparpenalty\predisplaypenalty
     \@endparpenalty\postdisplaypenalty
     \refstepcounter{equation}%
     \trivlist \item[]\leavevmode
       \hb@xt@\linewidth\bgroup $\m@th
         \displaystyle
         \hskip\mymathindent}%
        {$\hfil 
         \displaywidth\linewidth\hbox{\@eqnnum}%
       \egroup
     \endtrivlist}

\newcommand{\dext}{\mathrm{d}}

\newcommand{\del}{\partial}

\numberwithin{equation}{section}

\theoremstyle{plain}
\newtheorem{thm}{Theorem}[section]
\newtheorem{prp}[thm]{Proposition}
\newtheorem{lem}[thm]{Lemma}

\newtheorem{dfn}[thm]{Definition}
\newtheorem{exa}[thm]{Example}

\newtheorem{rmk}[thm]{Remark}

\newcounter{aequation} 

\author{Max Reinhold Jahnke
\footnote{
     I would like to thank Sönke Rollenske and Oliver Goertsches for their helpful discussions on spectral sequences and homogeneous spaces, respectively. This work was supported by the Deutsche Forschungsgemeinschaft (DFG), grant JA 3453/1-1.
}
\\ \href{mailto:jahnkem@mathematik.uni-marburg.de}{jahnkem@mathematik.uni-marburg.de} \\
\\
Fachbereich Mathematik und Informatik, \\
Philipps-Universität Marburg, \\
Marburg, Hessen, Germany.
}

\title{Closed elliptic structures on compact semisimple Lie groups}

\begin{document}

    \maketitle

    \begin{abstract} 
        In this work, we prove that, under a topological condition, the cohomology associated with left-invariant elliptic structures on compact semisimple Lie groups can be computed using only left-invariant forms. This reduces the analytical problem to a purely algebraic one, while also providing a generalization of the classic works of Chevalley and Eilenberg \cite{chevalley1948cohomology} on the de Rham cohomology of compact Lie groups and of Pittie \cite{pittie1988dolbeault} on the Dolbeault cohomology of compact semisimple Lie groups to the context of elliptic structures. We use spectral sequences as our primary tool, which facilitates the construction of an isomorphism between the left-invariant differential complex and the usual differential complex.
    \end{abstract}



    \section{Introduction}

        A fundamental problem in the theory of partial differential equations is the determination of the solvability of linear systems of vector fields. The theory of involutive structures provides us with a way to translate such solvability problems into the computation of cohomology spaces associated with complexes of differential operators. This approach allows us to treat in a unified way many interesting and important examples of differential complexes, such as the de Rham, Dolbeault, and tangential Cauchy-Riemann complexes. It also emphasizes the geometric and topological properties of the underlying manifold, thus facilitating the application of tools from other mathematical fields, such as algebraic topology and functional analysis. 
        
        When the underlying manifold is a compact Lie group, it is possible to connect the analytical problems related to the solvability of certain families of left-invariant vector fields with the algebraic properties of the Lie algebra of the manifold. By using this approach, we can transform the analytical problem of computing the cohomology spaces into simple problems in linear algebra. This approach has been done in several cases such as the de Rham \cite{chevalley1948cohomology}, left-invariant complex structures on compact Lie groups \cite{pittie1988dolbeault}, Levi-flat CR structures \cite{jacobowitz2021cohomology} and elliptic structures with certain topological and algebraic conditions \cite{jahnke2021cohomology}. In particular, left-invariant involutive structures have received increasing attention, as shown by the recent contributions of Charbonnel and Ounaïes-Khalgui \cite{charbonnel2004classification}, Korman \cite{korman2014elliptic}, Araújo \cite{araujo2019global}, and Bor and Jacobowitz \cite{bor2021complex}.
        
        
        %

        The following theorem is the main result of this paper.
        It significantly improves the main result from \cite{jahnke2021cohomology} by adapting the topological condition and removing the algebraic one.
        Additionally, it ensures a complete understanding of the cohomology classes associated with left-invariant elliptic structures in terms of the properties of the corresponding Lie algebras.
        
        \begin{thm}
        \label{thm:main}
            Suppose $G$ is a connected, compact, and semisimple Lie group, and let $\mathfrak v \subset \mathfrak g_\C$ be a left-invariant elliptic structure. Let $G_\C$ be the universal complexification of $G$. If $V = \exp_{G_\C}(\mathfrak v) \subset G_\C$ is closed, then every cohomology class in $H^{p,q}(G; \mathfrak v)$ has a left-invariant representative. Moreover, the inclusion of the left-invariant complex into the usual one induces an isomorphism in cohomology:
            \begin{equation}
                \label{eq:the_homomorphism}
                \phi :  H^{p,q} (\mathfrak g_\C; \mathfrak v) \to H^{p,q} (G; \mathfrak v).
            \end{equation}
        \end{thm}

        In Section \ref{sec:involutive_structures}, we establish notation and provide a brief overview of the theory of involutive structures. For a comprehensive introduction, we recommend consulting \cite[Chapters I, VII]{berhanu2008introduction} and \cite[Chapters I, III]{treves1992hypo}.
        In Section \ref{sec:basicdef}, we describe the connection between left-invariant involutive structures and Lie algebras, and the precise definition of the map $\phi$.

        

        The compactness of the Lie group is necessary for our context because we need the cohomology spaces $H^{p,q}(G; \mathfrak v)$ to be finite-dimensional. 
        If the Lie group is not compact, these cohomology spaces may not be finite-dimensional. 
        The semi-simplicity condition is needed as it allows us to apply a result due to Bott that requires that a certain homogeneous manifold is simply connected, and compact. This result plays an important role in our analysis, and it requires a certain level of structure on the Lie group and its Lie algebra. Finally, by assuming that $V$ is a closed subset of $G_\C$, we can take full advantage of Bott's theorem and its implications for our analysis.

        Theorem \ref{thm:main} covers the case $\mathfrak v = \mathfrak g_\C$ and the case where $\mathfrak v$ defines a left-invariant complex structure on $G$ \cite{pittie1988dolbeault}. The following class of examples shows the existence of a family of elliptic structures that satisfies the assumptions of Theorem \ref{thm:main}.
        
        \begin{exa}
            Let $G$ be a connected, compact, and semisimple Lie group, and let $\mathfrak t \subset \mathfrak g$ be a maximal abelian subalgebra. We denote by $\Delta$ the set of roots of $\mathfrak t_\C$ in $\mathfrak g_\C$, and by $\Delta_+$ a maximal subset of positive roots of $\Delta$. For each $\alpha \in \Delta$, let $\mathfrak g_\alpha$ be the corresponding eigenspace. It is easy to show that $\mathfrak v = \mathfrak t_\C \oplus \bigoplus_{\alpha \in \Delta_+} \mathfrak g_\alpha$ is an elliptic subalgebra of $\mathfrak g_\C$. We consider the universal complexification $G_\C$ of $G$ and note that $V \doteq \exp_{G_\C}(\mathfrak v) \subset G_\C$ is closed \cite[Claim 23.45, Proposition 26.4]{fulton1991representation}.
        \end{exa}
        
        Note that the example above can easily be adapted for a family of parabolic Lie algebras.


        The proof of Theorem \ref{thm:main} is divided into four steps. The first step of the proof is done in Section \ref{sec:involutive_structures}, where we briefly review the basic concepts of left-invariant involutive structures and set up the proof of Theorem \ref{thm:main} in terms of the surjectivity of a homomorphism. 

        The next step is to establish a connection between Lie algebra cohomology and a spectral sequence. In Section \ref{sec:cohomology_of_lie_algebras}, we provide a brief overview of Lie algebra cohomology and explain how it is related to the cohomology spaces $H^{p,q}(G; \mathfrak v)$ associated with the left-invariant involutive structures. We aim to prepare for the use of the Hochschild-Serre spectral sequence, which allows us to compute the second page of the spectral sequence associated with the Lie algebra cohomology induced by $\mathfrak v$.

        Then, by using an essentially real subalgebra of $\mathfrak v$, we construct an auxiliary Lie group and a homogeneous space. This is done in Section \ref{sec:leray_spectral_sequence}, where we explain how the Leray spectral sequence, which is a powerful tool in algebraic topology, relates the cohomology of the differential complex associated with the elliptic structure to the cohomology of the auxiliary group and the cohomology of the homogeneous space. 

        Finally, the fourth and last step of the proof is done in Section \ref{sec:the_proof}, where we show how to use results from Chevalley-Eilenberg \cite{chevalley1948cohomology} and Bott \cite{bott1957homogeneous} to combine the Hochschild-Serre and Leray spectral sequences. Essentially, the idea consists in constructing an isomorphism between the second pages of both spectral sequences and use it to conclude that there is an isomorphism between the limit terms.
        
        We conclude the paper with Section \ref{sec:open_problems}, where we discuss several open problems that arise from our work and suggest directions for future research. In particular, we highlight a few natural extensions of our results, including the study of cohomology classes associated with elliptic and non-elliptic Lie subalgebras and the study of left-invariant structures on homogeneous manifolds.

    
    \section{Involutive structures}
    \label{sec:involutive_structures}

        In this section, we present an overview of the essential ideas concerning involutive structures, their differential complexes, and the corresponding cohomology spaces. This serves both to explain the notation used throughout the paper and to improve its readability. The primary sources for this section are Chapters I and VII of ``Introduction to Involutive Structures'' by Berhanu, Cordaro, and Hounie \cite{berhanu2008introduction}.
        
        Let $\Omega$ be a smooth, oriented, and connected manifold of dimension $N \geq 2$. An involutive structure over $\Omega$ is a smooth subbundle $\mathcal V$ of $T_{\C} \Omega$ that satisfies the Frobenius condition: $[\mathcal V, \mathcal V] \subset \mathcal V$. Consider a finite-dimensional vector bundle $E$ over $\Omega$ and an open set $U \subset \Omega$. The set of smooth sections of $E$ over $U$ is denoted by $C^\infty(U; E)$.
        
        There is a natural differential complex associated with every involutive structure $\mathcal{V}$. For the convenience of the reader, we briefly recall the construction of such complexes. For each $p, q \in \Z_+$, we denote by $T^{p,q}_x$ the subspace of $(\Lambda^{p+q} T_\C^* \Omega)_x $ consisting of linear combinations of exterior products $u_1 \wedge \ldots \wedge u_{p+q}$ with $u_j \in (T_\C^*\Omega)_x$ for $j = 1, \ldots, p+q$ and at least $p$ of these factors belong to $T'_x$. Note that $T^{p+1,q-1}_x \subset T^{p,q}_x$ and we can define $\Lambda^{p,q}_x \doteq  T^{p,q}_x / T^{p+1,q-1}_x$ and $\Lambda^{p,q} = \bigcup_{x \in \Omega} \Lambda^{p,q}_x $ is a smooth vector bundle over $\Omega$.
        
        Since $\mathcal{V}$ is involutive we can prove that the exterior derivative maps smooth sections $T^{p,q}$ into smooth sections of $T^{p,q+1}$. Therefore, there exists a unique operator $\dext'_{p,q}$ which makes the following diagram commutative:
        $$
            \xymatrix{C^\infty(\Omega; T^{p,q}) \ar[r]^{\dext} \ar[d]_{\pi}& C^\infty(\Omega; T^{p,q+1}) \ar[d]^{\pi} \\
            C^\infty(\Omega; \Lambda^{p,q}) \ar[r]_{\dext'_{p,q}} & C^\infty(\Omega; \Lambda^{p,q+1}).}
        $$
        
        To simplify the notation, we denote by $\dext'$ the operator $\dext'_{p,q}$, or $\dext'_{\mathcal V}$ when it is necessary to emphasize the associated involutive structure. Note that the operator $\dext'$ maps smooth sections of $\Lambda^{p,q}$ into smooth sections of $\Lambda^{p,q+1}$, and it also holds that $\dext' \circ \dext' = 0$. Therefore, for each $p \in \{0, \ldots, m\}$, the operator $\dext'$ defines a complex of $\mathbb{C}$-linear mappings
        \begin{equation}
        \label{eq:h_complex}
           (\mathcal C^\infty(\Omega; \Lambda^{p,*}), \dext')
        \end{equation}
        and cohomology spaces denoted by
        \begin{equation}
        \label{eq:cohomology}
           H^{p,*}(\Omega; \mathcal V).
        \end{equation}
        
        \subsection{Left-invariant structures}
        \label{sec:basicdef}

        In this section, we consider a special class of involutive structures on Lie groups. As a reference for the basic concepts of Lie groups, we rely on the book ``Structure and Geometry of Lie Groups'' by Hilgert and Neeb \cite{hilgert2011structure}.
        
        Let $G$ be a Lie group with a Lie algebra $\mathfrak{g}$ and let $L_x : G \to G$ be the left-multiplication on $G$ by $x$. Recall that a vector bundle $\mathcal V \subset \mathbb{C} T G$ is called \emph{left-invariant} if $(L_x)_*X_g \in \mathcal V_{xg}, X_g \in \mathcal V_g.$
    
        If $\mathcal V$ is a left-invariant involutive vector bundle, there exists a corresponding Lie algebra $\mathfrak v \subset \mathfrak g_\C$ defined by the restriction of $\mathcal V$ to the identity of $G$. Note that this is a one-to-one correspondence, i.e., given a Lie subalgebra $\mathfrak v \subset \mathfrak g_\C$, we define by left-translation an involutive vector bundle $\mathcal V \subset T_\C G$.
        In this paper, we will always represent left-invariant involutive structures by their corresponding Lie algebras.
        Note that the rank of the vector bundle $\mathcal V$ is just the dimension of the Lie algebra $\mathfrak v$ and that the commutator bracket and the Lie algebra bracket agree.
    
        We are interested in the algebraic properties of the space $H^{p,q}(G; \mathfrak v)$. The properties of this space are influenced both by the inclusion of the Lie algebra $\mathfrak v$ into the algebra $\mathfrak g_\C$ and by the intrinsic algebraic properties of $\mathfrak v$, as well as by the topological aspects of the compact Lie group $G$. To highlight the algebraic properties of the locally integrable structures on $G$, we use the language of Lie algebras to characterize important classes of involutive structures relevant to this work:
        \begin{dfn}
            We say that a Lie subalgebra $\mathfrak v \subset \mathfrak g_\C$ defines: an \emph{elliptic structure} if $\mathfrak v + \bar {\mathfrak v} = \mathfrak g_\C$, and in this case, $\mathfrak v$ is called an \emph{elliptic algebra}; a \emph{complex structure} if $\mathfrak v \oplus \bar {\mathfrak v} = \mathfrak g_\C$, and in this case, $\mathfrak v$ is called a \emph{complex algebra}; and an \emph{essentially real structure} if $\mathfrak v = \bar {\mathfrak v}$, and in this case, $\mathfrak v$ is called a \emph{essentially real algebra}.
        \end{dfn}

        For example, if $\mathfrak v = \mathfrak g_\C$, the complex in \eqref{eq:h_complex} becomes the de Rham complex. In this case, we have an elliptic structure, which also happens to be an essentially real structure, and the operator $\dext'$ is the usual exterior derivative $\dext$. Now, if $\mathfrak v \oplus \overline{\mathfrak v} = \mathfrak g_\C$, we have a complex structure that is also elliptic. In this case, the operator $\dext'$ is the operator $\bar \del$, and the complex in  \eqref{eq:h_complex} is then called the Dolbeault complex.

        Since the exterior derivative of a left-invariant differential form is left-invariant, the operator $\dext'$ also maps left-invariant sections of $\Lambda^{p,q}$ to left-invariant sections of $\Lambda^{p,q+1}$. Thus, for each $p \in \{0, 1, \ldots, m\}$, we have another complex
        \begin{equation}
            \label{eq:h_complex_left_inv}
                (\mathcal C^\infty_L(G; \Lambda^{p,*}), \dext') 
        \end{equation}
        with cohomology spaces denoted by $H^{q}_L (G; \mathfrak v)$. In the next section, we briefly review the essential definitions and results of Lie algebra cohomology, which are crucial for constructing the cohomology spaces $H^{q}_L (G; \mathfrak v)$ using only the Lie algebras $\mathfrak{g}_\C$ and $\mathfrak v.$
        
        Note that, since $\mathcal C^\infty_L(G; \Lambda^{p,q}) \subset \mathcal C^\infty(G; \Lambda^{p,q})$ for all $p, q$, we have a natural homomorphism
        \begin{equation*}
            \phi :  H^{p,q} (\mathfrak g_\C; \mathfrak v) \to H^{p,q} (G; \mathfrak v).
        \end{equation*}
        
        Assuming that $G$ is compact, an averaging technique can be used to establish the injectivity of the homomorphism $\phi$. A detailed proof of this fact can be found in \cite[Lemma 5.2]{jahnke2021cohomology}.
        
        Let $G_\C$ be the universal complexification of $G$. Our current goal is to prove the following proposition:

        \begin{prp}
        \label{prp:simply_connected_case}
            Let $G$ be a connected, simply connected, and compact Lie group. Suppose that $\mathfrak v$ is an elliptic subalgebra of $\mathfrak g_\C$ and that the set $\exp_{G_\C}(\mathfrak v)$ is closed in $G_\C$. Then the homomorphism \eqref{eq:the_homomorphism} is surjective.
        \end{prp}
    
        The previous proposition implies Theorem \ref{thm:main} in the case where $G$ is a connected, simply connected, and compact Lie group. \todo{GA: Não entendi, esta é exatamente a hipótese da Proposição 2.2} The general case follows immediately since the left-invariant cohomology is preserved under passage to finite coverings. For details on the last claim, the reader can check \cite[Proposition 4.10 and Remark 4.11]{jahnke2021cohomology}.
          
    
    \section{Lie algebras cohomology}
    \label{sec:cohomology_of_lie_algebras}
            
        In this section, we provide a brief overview of the fundamental definitions and results of Lie algebra cohomology. For a comprehensive introduction to Lie algebra cohomology, we refer the interested reader to \cite[Chapter 7.5]{hilgert2011structure}.

        Let $\mathfrak{g}$ be an arbitrary Lie algebra and let $M$ be a $\mathfrak{g}$-module, i.e., a vector space equipped with a Lie algebra homomorphism $\rho: \mathfrak{g} \to \mathfrak{gl}(M)$. For $X \in \mathfrak{g}$ and $x \in M$, we denote $\rho(X)x$ by $X \cdot x$. We define $M^{\mathfrak{g}}$ to be the subspace of $M$ consisting of all $x \in M$ for which $X \cdot x = 0$ for all $X \in \mathfrak{g}$. Let $C^p(\mathfrak{g}; M)$ denote the set of all alternating multilinear $p$-forms on $\mathfrak{g}$ with values in $M$, and let $C^*(\mathfrak{g}; M) = \bigoplus_{p \geq 0} C^p(\mathfrak{g}; M)$ be the direct sum of all $C^p(\mathfrak{g}; M)$. We identify $C^0(\mathfrak{g}; M)$ with $M$ and we note that each $C^p(\mathfrak{g}; M)$ is a vector space. 

        For any non-negative integer $p$, let $u$ be an alternating $p$-linear form on $\mathfrak{g}$ with values in $M$, and let $X_1, \ldots, X_{p+1}$ be elements of $\mathfrak{g}$. We define the operator $\dext: C^p(\mathfrak{g}; M) \to C^{p+1}(\mathfrak{g}; M)$ by the formula:
        \begin{equation}
        \label{eq:algebraic_diff_op}
            \begin{aligned}
                 \dext u(X_1, \ldots, X_{p+1}) & = \sum_{j=1}^{p+1} (-1)^{j+1} X_j \cdot u(X_1, \ldots, \hat{X}j, \ldots, X_{p+1}) \\
                    & + \sum_{j < k} (-1)^{j+k+1} u([X_j, X_k], X_1, \ldots, \hat X_j, \ldots, \hat X_k, \ldots, X_p),
            \end{aligned}
        \end{equation}
        in which $\hat{X}_j$ denotes the omission of $X_j$.
            
        It can be shown that the operator $\dext$ satisfies $\dext \circ \dext = 0$, and so we have a cochain complex $(C^*(\mathfrak{g}; M), \dext)$. The cohomology spaces of this complex are denoted by $H^*(\mathfrak{g}; M)$.

        If $\mathfrak{v}$ is an ideal of $\mathfrak{g}$ (including the case $\mathfrak{v} = \mathfrak{g}$), we can define a $\mathfrak{g}$-module structure on the spaces $C^p(\mathfrak{v}; M)$. In fact, $C^0(\mathfrak{v}; M) = M$ already has a $\mathfrak{g}$-module, and for $p > 0$, $u \in C^p(\mathfrak{v}; M)$, $X \in \mathfrak{g}$, and $Y_1, \ldots, Y_p \in \mathfrak{v}$, we define the $\mathfrak{g}$-action on $u$ by $$(\mathcal{L}_X u)(Y_1, \ldots, Y_p) = X \cdot u(Y_1, \ldots, Y_p) - \sum_{i=1}^p u(Y_1, \ldots, [X,Y_i], \ldots, Y_p).$$
                        
        We have $\mathcal L_Y u \in C^p(\mathfrak v; M)$ for all $u \in C^p(\mathfrak v; M)$ and that $\mathcal L$ depends linearly on $Y$ and on $u$.
            
        \subsection{The cohomology induced by subalgebras}
            
            In this section, we review some basic concepts of Lie algebra cohomology induced by a subalgebra. We note that there are different  definitions of Lie algebra cohomology with respect to subalgebras. In this paper, we will make use of \emph{two} different notions. First, we introduce the notion corresponding to the algebraic version of the cohomology associated with left-invariant involutive structures.
        
            Let $\mathfrak v$ be a subalgebra of $\mathfrak g$. To define the cohomology induced by $\mathfrak v$, for $q \geq 0$, we define $N^{0,q}_{\mathfrak v}(\mathfrak g; M) = C^q(\mathfrak g; M)$ and for $p \geq 0$, we define $N^{p,q}_{\mathfrak v}(\mathfrak g; M)$ as the set of all $u \in C^{p+q}(\mathfrak g; M)$ that vanish when evaluated on $q+1$ elements of $\mathfrak v$. By convention, if $q < 0$, we set $N^{p,q}_{\mathfrak v}(\mathfrak g; M) = {0}$.
        
            Given the inclusion 
            \begin{equation}
            \label{eq:relative}
                N^{p,q}_{\mathfrak v}(\mathfrak g; M) \subset N^{p+1,q-1}_{\mathfrak v}(\mathfrak g; M),
            \end{equation}
            we can define the quotient space
            $ C^{p,q}_{\mathfrak v}(\mathfrak g; M) \doteq N^{p,q}_{\mathfrak v}(\mathfrak g; M) / N^{p+1,q-1}_{\mathfrak v}(\mathfrak g; M)$. Furthermore, we have that
            \begin{equation}
            \label{eq:closed_by_dext}
                \dext N^{p,q}_{\mathfrak v}(\mathfrak g; M) \subset N^{p,q+1}_{\mathfrak v}(\mathfrak g; M),
            \end{equation} the operator $\dext$ induces the operator
            $$ \dext'_{\mathfrak v} : C^{p,q}_{\mathfrak v}(\mathfrak g; M) \to C^{p,q+1}_{\mathfrak v}(\mathfrak g; M).$$
        
            Thus, for each $p$, we obtain a cochain complex $(C^{p,q}_{\mathfrak v}(\mathfrak g; M), \dext')$ whose cohomology spaces are denoted by $H^{p,*}_{\mathfrak v}(\mathfrak g; M)$.
            
            \subsection{The cohomology relative to a subalgebra}
        
            In addition to the cohomology associated with left-invariant involutive structures, we introduce a second notion of cohomology. This notion is essential for the application of the Hochschild-Serre spectral sequence.
            
            Let $\mathfrak v$ be a subalgebra of $\mathfrak g$. We define $C^p(\mathfrak g, \mathfrak v; M) \doteq N^{0,q}_\mathfrak v(\mathfrak g; M)$ and we note that by \eqref{eq:closed_by_dext} and \eqref{eq:relative} we have that $\dext C^p(\mathfrak g, \mathfrak v; M) \subset C^{p+1}(\mathfrak g, \mathfrak v; M)$. This leads us to define another complex  $(C^*(\mathfrak g, \mathfrak v; M), d)$ with cohomology spaces denoted by $H^*(\mathfrak g, \mathfrak v; M)$, which is called the cohomology of $\mathfrak g$ relative to $\mathfrak v$ with values on $M$.
            
            We define the space $C^p(\mathfrak g / \mathfrak v; M)$ to be the set of all alternating multilinear $p$-forms on the vector space $\mathfrak g / \mathfrak v$ with values in $M$. This space has a structure of a $\mathfrak v$-module. In fact, if $X \operatorname{mod} \mathfrak v \in \mathfrak g / \mathfrak v$ and $Y \in \mathfrak v$, the bracket $[X  \operatorname{mod} \mathfrak v,Y] = [X,Y] \operatorname{mod} \mathfrak v$ is well defined and so, for $X \in \mathfrak v$ and $Y_1 \operatorname{mod} \mathfrak v, \ldots, Y_p \operatorname{mod} \mathfrak v \in \mathfrak g / \mathfrak v$, we define
            $$ \begin{aligned}
                (\mathcal L_X u)(Y_1 \operatorname{mod} \mathfrak v, \ldots  , Y_p \operatorname{mod} \mathfrak v ) &= X \cdot u(Y_1 \operatorname{mod} \mathfrak v, \ldots, Y_p \operatorname{mod} \mathfrak v) \\
                &- \sum_{i=1}^n u(Y_1 \operatorname{mod} \mathfrak v, \ldots, [X,Y_i] \operatorname{mod} \mathfrak v, \ldots, Y_p \operatorname{mod} \mathfrak v).
            \end{aligned}$$
            
            The following theorem is used later in this paper. For a proof, see \cite{hochschild1953cohomology}.
            
            \begin{thm}
                 The complexes $(C^{p,q}_{\mathfrak v}(\mathfrak g; M), \dext_{\mathfrak v})$ and $(C^q(\mathfrak v; C^p(\mathfrak g / \mathfrak v; M)), \dext)$ are isomorphic and we have $$ H^{p,q}_{\mathfrak v}(\mathfrak g; M) = H^{q}(\mathfrak v; C^p(\mathfrak g / \mathfrak v ; M)). $$
            \end{thm}
        
            \subsection{An Application of the Hochschild-Serre Spectral Sequence}
        
            In order to establish certain algebraic properties of the cohomology spaces, we will use spectral sequences. It is outside the scope of this paper to review spectral sequences. We refer interested readers to \cite{bott982differential} for an excellent introduction to the topic.
            In the following, we describe in detail how to use spectral sequences to obtain the results presented here.
            
            Let $0 \leq q \leq \dim_\C \mathfrak g$ and $0 \leq p \leq q+1$. Let $F^p C^q(\mathfrak g; M) = N^{p, q-p}(\mathfrak g; M)$. We have that $F^p C^*(\mathfrak g; M)$ is a graded subspace of $C^*(\mathfrak g; M)$, and we have a descending sequence of spaces
            $$C^q(\mathfrak g; M) = F^0 C^q(\mathfrak g; M) \supset \cdots \supset F^{q+1} C^q(\mathfrak g; M) = 0.$$ From this filtration, we obtain a spectral sequence $(E^{*,*}_r, \dext_r)$ with limit term associated with $H^*(\mathfrak g; M)$.
            
            We now consider the case where $\mathfrak g$ is the Lie algebra of a compact Lie group. Let $\mathfrak v \subset \mathfrak g_\C$ be an elliptic subalgebra, and let $\mathfrak k = \mathfrak v \cap \overline{\mathfrak v}$.
            
            Since $\mathfrak g$ is the Lie algebra of a compact Lie group, it is equipped with an ad-invariant inner product $\langle,\rangle$. We can extend this inner product to a Hermitian inner product on $\mathfrak g_\C$, which now satisfies the relation:
        
            \begin{equation}
            \label{eq:not_ad_inv}
                \langle [X,Y], Z \rangle = - \langle Y, [\overline X,Z] \rangle, \quad \forall X,Y,Z \in \mathfrak g_\C.
            \end{equation} 
        
            Note that this extension is not ad-invariant.

            \todo{GA: Omitir a demonstração?}
            \begin{lem}
                The adjoint representation of $\mathfrak k$ in $\mathfrak v$ is completely reductive, i.e., every $\mathfrak{k}$-invariant subspace of $\mathfrak{v}$ has a $\mathfrak{k}$-invariant complement.
            \end{lem}
            
            \begin{proof}
                To show that $\mathfrak v$ is a completely reductive $\mathfrak k$-module, we need to prove that every $\mathfrak k$-submodule $M$ of $\mathfrak v$ has a $\mathfrak k$-invariant complement. Let $M \subset \mathfrak v$ be $\mathfrak k$-submodule of $\mathfrak v$ and let $M^\perp$ be the orthogonal complement of $M$ with respect to the hermitian inner product \eqref{eq:not_ad_inv}, that is, $M^\perp = \{ Z \in \mathfrak v; \langle Z, X \rangle = 0 \text{ for all } X \in K\}.$
                
                By construction, $\mathfrak v$ can be written as a direct sum of vector spaces $\mathfrak v = M \oplus M^\perp$. Now we need to show that $M^\perp$ is $\mathfrak k$-invariant. Let $X \in \mathfrak k$ and $Z \in M^\perp$. Note that we can assume that $X$ is a real vector field, and as a result we have that $\overline{X} = X$. Therefore we have that $\langle [X, Z], Y \rangle = -  \langle Z, [X,Y] \rangle = 0$ because $[X,Y] \in \mathfrak k$ and we obtain $[X,Z] \in M^\perp$.
            \end{proof}
            
            Now, consider the complex $(C^{p,*}_{\mathfrak v}(\mathfrak g), \dext') $ which is isomorphic to the complex $(C^*(\mathfrak v; C^p(\mathfrak g / \mathfrak v)), \dext)$. By considering the filtration $F^* C^*(\mathfrak v; C^p(\mathfrak g/ \mathfrak v; \C))$, we obtain an spectral sequence $(E^{*,*}_r, \dext_r)$ which, by Theorem 11 of \cite{hochschild1953cohomology}, has the following $E_2$ term
            \begin{equation}
            \label{eq:hs_page_2}
                E^{j,k}_2  \cong H^j(\mathfrak k; \C) \oplus H^k(\mathfrak v, \mathfrak k, C^p(\mathfrak g/ \mathfrak v; \C)).
            \end{equation}


    
    \section{Leray spectral sequence}
    \label{sec:leray_spectral_sequence}

    In this section, we establish the necessary foundations for the application of the Leray spectral sequence by introducing a homogeneous space associated with the Lie group $G$. This homogeneous space is used in two ways. First, we show how the natural complex structure of the homogeneous space and the differential structure of a subgroup allow us to understand the local structure of the elliptic structure on $G$, which is important for understanding the cohomology of certain sheaves that appear in the Leray spectral sequence. Then, we compute the second page of the Leray spectral sequence and connect it to the Hochschild-Serre spectral sequence. 
    


    Let $G$ be a compact connected Lie group and let $\mathfrak v \subset \mathfrak{g}_\C$ be an elliptic Lie algebra. We consider the universal complexification $G_\C$ of $G$, and we define the group $H = \exp{G_\C}(\mathfrak v)$, which we assume to be a closed subgroup of $G_\C$. Then, we construct the quotient space $\Omega = G_\C / H$ with the canonical projection map $\pi: G_\C \to \Omega$. Since $G$ is compact, we can identify $G$ with a compact subgroup of its universal complexification $G_\C$. Note that we can take $K \doteq H \cap G$ and we have that $K$ is a closed subgroup of $G$ with $\mathfrak k = \mathfrak v \cap \mathfrak g$ as its Lie algebra.

    Next, we show that $\Omega$ can also be expressed as $G/K$. For this, we consider the homogeneous space $\Omega' = G/K$, which is equipped with the projection map $\phi : G \to \Omega'$. It suffices to show that $\Omega'$ is diffeomorphic to $\Omega$. 

    We define a map $F:\Omega' \to \Omega$ as $gK \mapsto gH$. The map $F$ is clearly well-defined and a injective immersion. From the fact that $\Omega$ is compact, we have that $F(\Omega')$ is an embedded manifold. To show that the map is surjective, it suffices to show that the dimensions of $\Omega'$ and $\Omega$ are equal. Since $\mathfrak v$ is elliptic, we have $\mathfrak g_\C = \mathfrak v + \overline{\mathfrak v}$, implying that $$\dim_\R \mathfrak g = \dim_\C \mathfrak g_\C = \dim_\C \mathfrak k_\C + 2 \dim \mathfrak v / \mathfrak k_\C \quad$$  and $$\quad \dim_\R \mathfrak v = \dim_\R \mathfrak k_\C - \dim_\R \mathfrak v / \mathfrak k_\C,$$
    from which we easily conclude that $\dim \Omega = \dim \Omega'.$

    Therefore, $F(\Omega')$ is open, and since $\Omega$ is connected and $F$ is injective and smooth, we conclude that $F(\Omega') = \Omega$ and $F$ is a diffeomorphism. We can thus identify $\Omega'$ with $\Omega$ via $F$. This means that the action of $G$ on $\Omega'$ can be lifted to a transitive action on $\Omega$, which now has two equivalent descriptions:
    \begin{equation}
    \label{eq:two_descriptions}
        \Omega = G_\C / H = G / K.
    \end{equation} 

    \begin{rmk}
        Note that if we additionally assume that $G$ is simply connected, then $G_\C$ is also connected and simply connected. Since $H$ is a connected subgroup of $G_\C$, the quotient space $G_\C/H$ is also connected. In fact, by the long exact sequence of homotopy groups, we obtain that $\pi_1(G_\C/H) \cong \pi_1(\Omega) \cong 0$, which means that $\Omega$ is simply connected.
    \end{rmk}

    Consider the sheaf $\mathcal S^p$ obtained from the presheaf $U \mapsto \mathcal S^p(U) = \{f \in C^\infty(U; \Lambda^{p,0}); \dext' f = 0\}$ for $U \subset G$ open. Since elliptic structures are locally solvable, we have that the sheaf cohomology $H^q(G; \mathcal S^p)$ and the usual cohomology spaces $H^q(G; \mathfrak v)$ are isomorphic. Therefore, we use sheaf cohomology to study the spaces $H^q(G; \mathfrak v)$. The idea now is to use the Leray spectral sequence as a tool to infer properties of $H^q(G; \mathfrak v)$ from the cohomology of $\Omega = G/K$, with $\Omega$ endowed with the natural complex cohomology induced by the quotient map $\phi : G \to \Omega = G/K$, and the cohomology of $K \subset G$ with the usual de Rham cohomology.

    Consider the sheaf $\mathcal H = \sum_{k=0}^\infty \mathcal H^k$ with $\mathcal H^k$ obtained from the presheaf
    \begin{equation}
    \label{eq:E_2_local}
        V \mapsto H^k(\phi^{-1}(V); \mathcal S^p|_{\phi^{-1}(V)}).
    \end{equation}
    The Leray spectral sequence is a spectral sequence $E_r$ with
    \begin{equation}
    \label{eq:E_2_term}
        E_2^{r,s} = H^r(\Omega; \mathcal H^s)
    \end{equation} with limit term associated with $H^*(G; \mathcal S^p)$. In the next sections, we will give an explicit description of $E_2^{r,s} = H^r(\Omega; \mathcal H^s)$.

    \subsection{The local structure of the elliptic structure}


    To provide an explicit description of the cohomology groups $E_2^{r,s} = H^r(\Omega; \mathcal H^s)$ in terms of $K$ and $\Omega$, we have to first understand the geometric structure of $G$ with respect to the elliptic structure $\mathfrak v$. To be precise, we prove that for all $x \in \Omega$, there exists a local section $\sigma$ for $\phi$ defined in a neighborhood of $x$, which is compatible with the complex structure of both $\Omega$ and $G$. Specifically, with $\Omega$ equipped with the complex structure $\mathcal W$ given by $\pi_*(\mathfrak v)$, the local sections $\sigma$ satisfy $\sigma_*(W) \subset \mathfrak v$. 

    \begin{lem}
    \label{lem:section_for_elliptic}
        Let $M$ and $N$ be connected and smooth manifolds. We assume that $M$ has an elliptic structure $\mathcal V$, and $N$ has a complex structure $\mathcal W$. Consider a map $f : M \to N$ such that $f_*(\mathcal{V}_x) = \mathcal{W}_{f(x)}$ for all $x \in M$. Then, for every $y \in N$, there exists an open neighborhood $U$ of $y$ and a map $\sigma : U \to M$ such that $f \circ \sigma(x) = x$ for all $x \in U$, and $\sigma_*(\mathcal W_x) \subset \mathcal V_{\sigma(x)}$ for all $x \in U$.
    \end{lem}

    \begin{proof}
        The local section is obtained by an application of the inverse function theorem. In order to make it explicit, we need suitable coordinates \cite[pg. 28]{berhanu2008introduction}. Let $\nu$ be the complex dimension of $N$ and let $m = \nu + n$ be the real dimension of $M$. Note that the hypothesis implies that $n > \nu$. Let $y \in N$, $x \in f^{-1}(y)$, and $V \subset M$ be an open neighborhood of $x$. We endow $V$ with coordinates $\phi = (z_1, \ldots, z_\nu, t_1, \ldots, t_{n - \nu}) : V \to V' \times V'' \subset \C^\nu \times \R^{n-\nu}$ with $V' \subset \C^\nu$ and $V'' \subset \R^{n - \nu}$ such that over $V$ the involutive structure $\mathcal V$ is generated by the vector fields $$\frac{ \del }{ \del \overline{z}_1}, \ldots, \frac{ \del }{ \del \overline{z}_\nu}, \frac{ \del }{ \del t_1}, \ldots, \frac{ \del }{ \del t_{n - \nu}}.$$
        
        Let $U \subset N$ be a neighborhood of $y$ with coordinates $ \psi = (w_1, \ldots, w_\nu) : U \to U' \subset \C^\nu$ such that over $U$ the involutive structure $\mathcal W$ is generated by the vector fields $$\frac{ \del }{ \del \overline{w}_1}, \ldots, \frac{ \del }{\del \overline{w}\nu}.$$
        
        Now, by using these coordinates, we write a local representation of $f$, namely $\tilde f: V' \times V'' \to U'$, as $\tilde f = \psi \circ f \circ \phi^{-1}$. Note that the hypothesis on the pushforward by $f$ implies that $$\tilde f_*(\del / \del t_j) \in \operatorname{span}_{C^\infty(U')} \{\del / \del \overline{w}_1, \ldots, \del / \del \overline{w}_\nu\}$$ and note that since $\tilde f$ and $\del / \del t_l$ are real, we have that $$\tilde f_*(\del / \del t_j) = \overline{\tilde f_*( \del / \del t_j)} \in \overline{\operatorname{span}_{C^\infty(U')} \{\del / \del \overline{w}_1, \ldots, \del / \del \overline{w}_\nu\}}$$ and thus $\tilde f_*(\del / \del t_j) = 0$ for all $j$, that is, $\tilde f$ does not depend on the variables $t_j's$.
    
        Let $J \tilde f (z)$ denote the Jacobian matrix of $\tilde f$. Since $\tilde f_*$ is surjective, by shrinking the open sets $V,V',V'',U$ and $U'$ if necessary, and by rearranging indices, we can assume that the first $l \times l$ block of the Jacobian matrix of $\tilde f$ is invertible. We write $z' = (z_1, \ldots, z_l)$ and since $\tilde f$ does not depend on $t$, we define the function $F(z', z'') = (\tilde f(z',z'',t_0),z'')$. From the inverse function theorem, there exists $F^{-1}$ and this function is holomorphic. Let $F'(z,t) = (F(z),t)$,  and note that this function is invertible. We define $\phi' = F' \circ \phi$ and, with $\phi'$ as a new coordinate on $V$, we have that the local representation of $f$ is the projection on $z'$. Therefore, we can define a holomorphic section for $f$ by defining $\sigma' (z',z'',t) = (z', z''_0, t_0)$ and then defining $\sigma : U \to M$ as $\sigma = \phi^{-1} \circ \tilde \sigma \circ \phi$. Thus, we have that $f \circ \sigma = \operatorname{Id}_U$ and by construction, since $\tilde \sigma$, $\phi$ and $\psi$ are compatible with the structures $\mathcal V$ and $\mathcal W$, we have that $\sigma_*(w) \in V_{\sigma(x)}$ for all $w \in W_x$.
    \end{proof}

    \begin{prp}
    \label{prp:principal_bundle}
        Let $G$ be a compact Lie group endowed with a left-invariant elliptic structure $\mathfrak v$, and let $\mathfrak{k} = \mathfrak{v} \cap \mathfrak g$ be its real part which we integrate into a group $K = \exp_G(\mathfrak k)$. We assume $K$ to be closed, and we define $\Omega = G/K$ and endow $\Omega$ with the complex structure defined by $\pi_*(\mathfrak v)$ and denoted by $\mathcal W$. Let $\Omega \times K$ be endowed with the elliptic structure given by $\mathcal W \oplus \mathfrak k_\C$.
        Then, there exists a finite covering $\mathfrak U$ of $\Omega$ and local sections $\sigma = \sigma_U$ for $\pi$ defined over $U \in \mathfrak U$ such that
        \begin{equation} \label{eq:local_elliptic_diffeomorphism}
            \Phi : (u,k) \in U \times K \mapsto \sigma(u)k \in \sigma(U)K = \pi^{-1}(U)
        \end{equation}
        satisfies $(\Phi_*)_{(u,k)}(X \oplus Y) \in \mathfrak{v}_{\sigma(u)k}$.        
    \end{prp}

    \begin{proof}
        By Proposition 15.2.4 of \cite{hilgert2011structure} and Lemma \ref{lem:section_for_elliptic}, with $M = G$ and $N = \Omega (=G/K)$, for every small open set $U \subset \Omega$ and any $k \in K$, there exists a diffeomorphism
        \begin{equation}
            \Phi : (u,k) \in U \times K \mapsto \sigma(u)k \in \sigma(U)K = \pi^{-1}(U).
        \end{equation}
    
        Now we have to prove that the diffeomorphism we constructed satisfies $$(\Phi_*)_{(u,t)}(X \oplus Y) \in \mathfrak{v}_{\sigma(u)k}$$
        for all $X \oplus Y \in \pi_*(\mathfrak{v})|_u \oplus \mathfrak k_\C|_k$. Let $\Psi = \Phi^{-1}$. Note that $\Psi : \pi^{-1}(U) \to U \times K$ is given by $\Psi(g) = (\pi(g), \psi(g))$ with $\psi : \pi^{-1}(U) \to K$ being a smooth function such that its differential is surjective at every point. We have $\Psi_*(\mathfrak{v}) = (\pi_*(\mathfrak{v}), \phi_*(\mathfrak{v})) \subset \pi_*(\mathfrak h) \oplus \mathfrak{k}_\C$. Since $\Psi$ is a diffeomorphism and the rank of $\mathfrak{v}$ is the same as the rank of  $\pi_*(\mathfrak v) \oplus  \mathfrak{k}_\C$, we have that $\Psi_*(\mathfrak{v}) = \pi_*(\mathfrak v) \oplus  \mathfrak{k}_\C$.
    \end{proof}

    Now, we can finally compute the $E_2$ of the Leray spectral sequence associated with our problem. In the following, $\Omega$ is endowed with a complex structure induced by $\pi_*(\mathfrak v)$ and $\mathcal O^p$ denotes the sheaf of holomorphic $p$-forms on $\Omega$.

    \begin{prp}
        Let $G$ be a connected, compact, and semisimple Lie group, and let $\mathfrak v \subset \mathfrak g_\C$ be an elliptic algebra. Let $G_\C$ be the universal complexification of $G$. If $V = \exp_{G_\C}(\mathfrak v) \subset G_\C$ is closed, then we consider the two descriptions of $\Omega = G_\C / V = G / K$. Consider the Leray spectral sequence associated with the quotient map $ \phi : G \to \Omega $. Then, then the $E_2$ term \ref{eq:E_2_term} is
        \begin{equation}
            H^r(\Omega; \mathcal H^s) \cong \sum_{r + s = q} H^r(\Omega; \mathcal O^p) \otimes H^s(K; \mathcal S_{\frk v}).
        \end{equation}
     \end{prp}

     \begin{proof}
        By Proposition \ref{prp:principal_bundle}, there exists a finite covering $\mathfrak U$ of $\Omega$ by open sets such that, for each $U \in \mathfrak U$, there is a smooth section $\sigma$ of the quotient map $\phi : G \to \Omega$ such that the following map $$ (u, k) \in U \times K \xmapsto{ \Phi } \sigma(u)k \in \sigma(U)K$$
        is a diffeomorphism preserving the involutive structures onto an open subset of $G$.
        Therefore, by using Künneth's formula, we have $$H^q(\phi^{-1}(U); \mathcal S^p|_{\phi^{-1}(U)}) = \sum_{r + s = q} H^r(U, \mathcal O^p) \otimes H^s(K,\mathcal S_{\frk v}).$$
        
        Thus,
        \begin{equation}
        \label{eq:leray_page_2}
            H^r(\Omega; \mathcal H^s) \cong \sum_{r + s = q} H^r(\Omega; \mathcal O^p) \otimes H^s(K; \mathcal S_{\frk v}).
        \end{equation}
     \end{proof}

        
    \section{The proof of Theorem \ref{thm:main}}
    \label{sec:the_proof}

    To prove Theorem \ref{thm:main}, it is sufficient to prove Proposition \ref{prp:simply_connected_case}. The latter can be shown by proving the existence of an isomorphism between page $E_2$ \eqref{eq:hs_page_2} of the Hochschild-Serre spectral sequence and page $E_2$ \eqref{eq:leray_page_2} of the Leray spectral sequence. Since such an isomorphism exists, it follows that the subsequent pages of the two sequences are also isomorphic \cite[Chapter 9, Section 1, Theorem 1]{spanier1994algebraic}, and thus the limit terms are isomorphic as well. This implies, by dimensional reasons, that the homomorphism \eqref{eq:the_homomorphism} is both injective and surjective.
    
    Now, we construct the isomorphism between page $E_2$ of the Hochschild-Serre spectral sequence and page $E_2$ of the Leray spectral sequence. By Theorems 2.2 and 2.4 of \cite{chevalley1948cohomology}, we have that $H^s(K,\mathcal S_{\frk v}) \cong H^s(\mathfrak k)$ and by Theorem II of \cite{bott1957homogeneous}, we have that $H^r(\Omega, \mathcal O^p) = H^r(\mathfrak v, \mathfrak k, C^p(\mathfrak g_\C / \mathfrak v)^*)$ with $v$ acting on $C^p(\mathfrak g_\C / \mathfrak v)^*$ via an adjoint action. The proof concludes by considering the product and sum of the isomorphisms.


    \section{Open problems}
    \label{sec:open_problems}
    \todo{GA: seria melhor omitir essa seção?}
    
     To encourage further research in this area, we present a list of open problems that naturally arise from our work. These include the study of cohomology classes associated with both elliptic and non-elliptic Lie subalgebras, as well as the study of left-invariant structures on homogeneous manifolds.

    \begin{enumerate}
              

        \item Is it possible to classify all elliptic Lie algebras $\mathfrak v$ such that the associated Lie group $V = \exp_{G_\C}(\mathfrak v)$ is closed?
        
        \item Note that Theorem \ref{thm:main} is valid for connected compact semisimple Lie groups endowed with closed elliptic structures. If we change $G$ for a torus, the same result is valid without assuming that the associated $V$ is closed (see \cite[Proposition 5.11]{jahnke2021cohomology}). Using the Künneth formula, a similar result can be obtained for connected compact Lie groups, provided that the elliptic structure is a product of elliptic structures. Is it possible to extend the main result to general closed elliptic structures on connected compact Lie groups without further algebraic restrictions on $\mathfrak v$?

        \item The question of whether every cohomology class in $H^{p,q}(G; \mathfrak v)$ has a left-invariant representative remains open for the general case, where $G$ is a connected, compact Lie group and $\mathfrak v \subset \mathfrak g_\C$ is an elliptic algebra and we do not impose any topological conditions.

        \item In our work, Bott's theorem was a crucial tool. This theorem provides an algebraic description of the Dolbeault cohomology for specific compact homogeneous manifolds. However, our study gave rise to two questions. Firstly, can we refine Bott's result by relaxing some of the topological constraints imposed on the manifold? Secondly, is it feasible to generalize Bott's theorem to encompass compact homogeneous manifolds with a left-invariant elliptic structure? 
        
    \end{enumerate}

    Answering these questions will require further investigation and may yield valuable insights into the interplay between the topological, analytical and algebraic structures on compact Lie groups and homogeneous manifolds.




    \bibliographystyle{alpha}
    \bibliography{bibliography}
    
\end{document}